\documentclass[12pt]{amsart}
\makeatletter
\DeclareOldFontCommand{\rm}{\normalfont\rmfamily}{\mathrm}
\DeclareOldFontCommand{\sf}{\normalfont\sffamily}{\mathsf}
\DeclareOldFontCommand{\tt}{\normalfont\ttfamily}{\mathtt}
\DeclareOldFontCommand{\bf}{\normalfont\bfseries}{\mathbf}
\DeclareOldFontCommand{\it}{\normalfont\itshape}{\mathit}
\DeclareOldFontCommand{\sl}{\normalfont\slshape}{\@nomath\sl}
\DeclareOldFontCommand{\sc}{\normalfont\scshape}{\@nomath\sc}
\makeatother

\usepackage[margin=1in]{geometry} 

\usepackage{times}
\usepackage[english]{babel}
\usepackage[utf8]{inputenc}
\usepackage{csquotes}
\usepackage{cancel}
\usepackage{hyperref}
\usepackage{amsmath,amssymb,amsthm,bbm}
\usepackage{cases}
\usepackage{amsfonts}
\usepackage{mathtools}
\usepackage{esint}
\usepackage{tabulary}
\usepackage{bm}
\usepackage{graphicx, color}
\usepackage{svg}
\usepackage{enumitem}
\usepackage{cancel}

\newcommand{\dmu}{\,{\rm d} \mu}

\newcommand{\dd}{\,{\rm d}} 

\newcommand{\dx}{{\rm d} x}

\newcommand{\dt}{{\rm d} t}
\newcommand{\ds}{{\rm d} s}
\newcommand{\dr}{{\rm d} r}

\newcommand{\dhx}[1]{{\rm d} \mathcal{H}^{#1}}

\newcommand{\dH}[1]{\,{\rm d} \mathcal{H}^{#1}}

\newcommand{\HH}{\mathcal{H}}

\newcommand{\V}{\mathbb{V}}
\newcommand{\vv}{\vec{v}}
\newcommand{\hh}{\vec{h}}
\newcommand{\m}{\mu_t}

\newcommand{\numberset}{\mathbb}
\newcommand{\R}{\numberset{R}}

\newcommand{\N}{\numberset{N}}

\newcommand{\Z}{\numberset{Z}}

\newcommand{\eps}{\varepsilon}

\newtheorem{theorem}{Theorem}[section]
\newtheorem{definition}{Definition}[section]
\newtheorem{proposition}{Proposition}[section]

\theoremstyle{definition}
\newtheorem{remark}{Remark}[section]

\title[The existence of Brakke-flow for volume preserving mean curvature flow]{Brakke inequality and the existence of Brakke-flow for volume preserving mean curvature flow}

\author[A. Chiesa] {Andrea Chiesa} 
\address[Andrea Chiesa]{University of Vienna, Faculty of Mathematics  and Vienna School of Mathematics,
                Oskar-Morgenstern-Platz 1, A-1090 Vienna, Austria}
\email{andrea.chiesa@univie.ac.at}
    
\author[K. Takasao] {Keisuke Takasao} 
\address[Keisuke Takasao]{ Department of Mathematics, Kyoto University, Kitashirakawa-Oiwakecho Sakyo
Kyoto 606-8502, Japan}
\email{k.takasao@math.kyoto-u.ac.jp}

	
\date{}

\begin{document}
\begin{abstract}
    In this paper, we propose a new notion of Brakke inequality for volume preserving mean curvature flow. We show the existence of integral varifolds solving the flow globally-in-time in the corresponding Brakke sense using the phase field method. Moreover such varifolds are solutions to volume preserving mean curvature flow in the $L^2$-flow sense as well. We thus extend a previous result by one of the authors \cite{TakasaoVPMCF}. 
\end{abstract}
\maketitle

\section{Introduction}

Let $d\in \N$, $d\geq 2$, and $\R^d$ be the $d$-dimensional Euclidean space. 
Assume that $T>0$ is a given final time and $U_t \subset \R^d$ is a bounded  open set with the smooth boundary 
$M_t:=\partial U_t$ for any $t \in [0,T)$.
The family of the hypersurfaces $\{ M_t \} _{t \in [0,T)}$ is called the \emph{volume preserving mean curvature flow}
if the normal velocity vector $\vec{v}$ satisfies
\begin{equation}
\vec{v}=\vec{h} - \left( \frac{1}{\mathcal{H}^{d-1} (M_t)} \int _{M_t} \vec{h} \cdot \vec{\nu} \,  \dH{d-1} \right) \vec{\nu},
\quad \text{on} \ M_t , \ t \in (0,T).
\label{vpmcf}
\end{equation}
Here, $\HH^{d-1}$ is the $(d-1)$-dimensional Hausdorff measure, and
$\vec{h}$ and $\vec{\nu}$ are the mean curvature vector and the inner unit normal vector of $M_t$, respectively.

The gradient flow equation \eqref{vpmcf} has been broadly studied in the literature: In the smooth setting, existence and convergence to a sphere as $t\rightarrow\infty$ have been established if the initial configuration $U_0$ is convex \cite{Gage,Huisken}. Moreover, \cite{EscherSimonett} showed short-time existence for smooth initial data $M_0$ and global-in-time existence whenever $M_0$ is close to a sphere in the sense of the little H{\"o}lder norm $h^{1+\alpha}$. On the other hand, a wide range of results is also present for weak notions of solution. In the context of the minimizing movements scheme, the sequence of time-discretized minimizers $\{U^k_t\}_{t\in [0,\infty)}$ converges to   a global weak solution $\{U_t\}_{t\in [0,\infty)}$ to \eqref{vpmcf} in the sense of the \emph{flat flow} \cite{MugnaiSeisSpadaro}, i.e., $\{U_t\}_{t\in [0,\infty)}$ is a sequence of Caccioppoli sets with constant volume $\mathcal{L}^{d}$ and monotone decreasing Hausdorff  measure of their (reduced) boundary $\HH^{d-1}(\partial^{*}U_t)$, satisfying the inequality $\mathcal{L}^d(U_s\triangle U_t)\leq c\sqrt{t-s}$ for every $0\leq s<t$, with $c=c(d,U_0)>0$. In \cite{Julin-Niinikoski}, the Alexandrov theorem and the asymptotic behavior of the flat flow were studied. Furthermore, if $d\leq 7$, \cite{MugnaiSeisSpadaro} also showed the existence of a solution to \eqref{vpmcf} on $[0,T]$ in the sense of distributions under the conditional convergence assumption 
\begin{equation}\label{convergence assumption}
    \lim_{k\rightarrow\infty}\int_0^T\HH^{d-1}(\partial^*U^k_t)\dt=\int_0^T\HH^{d-1}(\partial^*U_t)\dt
\end{equation}
introduced in \cite{LuckhausSturzenhecker}. Note that in general only $\liminf_{k\rightarrow\infty}\HH^{d-1}(\partial ^*U^k_t)\geq \HH^{d-1}(\partial ^*U_t)$. Analogous results under similar conditional convergence conditions have been established for distributional BV-solutions to \eqref{vpmcf} in the setting of the thresholding scheme \cite{LauxVolumePreserving}  and of the phase-field method \cite{LauxSimon}. The problem of the convergence of approximate solutions to weak solutions under assumption \eqref{convergence assumption} has been widely studied, see e.g. \cite{Almgren, LauxOtto,LuckhausSturzenhecker}. Recently, Laux proved the weak-strong uniqueness of the volume preserving mean curvature flow \cite{LauxWeakStrongUniqueness}. The existence of a viscosity solution to \eqref{vpmcf} under the assumption that the initial datu $U_0$ satisfies a suitable geometric condition, as well as uniform asymptotic convergence to a sphere, have been addressed in \cite{KimKwon}. Finally, in the setting of varifolds, \cite{TakasaoVPMCF} showed the existence of a solution to \eqref{vpmcf} in the sense of $L^2$-flow. 
Motivated by the fact that any canonical Brakke-flow is also an $L^2$-flow \cite[Thm. 4.3]{StuvardTonegawa}, the present paper aims at extending the result in \cite{TakasaoVPMCF} by showing that such solution satisfies a suitable Brakke inequality in the case of periodic boundary conditions.
Note that the Brakke flow
is an appropriate weak solution for analyzing the multi-phase mean curvature flow (without the volume constraint) because it is known for its existence and regularity properties (see \cite{Brakke, KasaiTonegawa,
StuvardTonegawa, Tonegawa2014, TonegawaBook} and the references therein). Moreover, recently \cite{StuvardTonegawa} showed the weak-strong uniqueness of the suitable weak solution, called the \emph{canonical Brakke flow}, by using \cite{Fischer}.
In light of the above, it is natural and meaningful to define and analyze the Brakke flow for \eqref{vpmcf} as well.

Let us start by recalling the notion of Brakke inequality and its motivation, we refer to \cite{TonegawaBook} for more details.  If $\{M_t\}_{t\in [0,T]}$ is a smooth $(d-1)$-dimensional hypersurface with normal velocity $\vec{v}$ and mean curvature $\vec{h}$, then it satisfies the following weighted first variation formula (cf. \cite[Remark 4.5]{Ecker} and \cite[Sec. 2.1]{TonegawaBook})
\begin{equation}\label{firstvariationformula}
    \frac{{\rm d}}{{\rm d} t}\int_{M_t}\phi\dH{d-1}=\int_{M_t}\left(\phi_t+(\nabla^\perp \phi - \phi\vec{h})\cdot \vec{v}\right)\dH{d-1},
\end{equation}
where $\nabla^\perp \phi$ is the component of $\nabla\phi$ normal to $M_t$. Notice that, since in \eqref{firstvariationformula} and in the following $\vec{v}$ is the normal velocity vector, we can simply write $\nabla \phi\cdot \vec{v}$.
The idea behind Brakke inequality relies on the identity \eqref{firstvariationformula} to define the normal velocity in a weak sense, by integrating in time and substituting the equality with an inequality. Namely, $\vec{v}$ is said to be (weak) normal velocity to the varifold $\{V_t\}_{t\in [0,T]}$ if   \begin{equation}\label{Brakke ineq}
            \|V_t\|(\phi(\cdot,t))\Big|_{t=t_1}^{t_2}\leq\int_{t_1}^{t_2}\int_{\R^d} \left(\phi_t+(\nabla\phi-\phi\hh)\cdot\vec{v}\right)\,{\rm d}\|V_t\|\dt
        \end{equation}
for every $\phi\in C^1_c({\R^d}{\times} [0,T];[0,\infty))$ and $0\leq t_1< t_2\leq T$. From the mathematical viewpoint, it is convenient to substitute the equality with an inequality, since \eqref{Brakke ineq} is in general obtained by approximation and the quantities converging to $\vec{v}$ and $\vec{h}$ do so weakly, so that one cannot expect to pass to the limit in the product $\vec{h}\cdot\vec{v}$ and preserve the identity. On the other hand, one can show that in the smooth setting, inequality \eqref{Brakke ineq} suffices to characterise the normal velocity. Indeed, let $\{M_t\}_{t\in[0,T]}$ be smooth hypersurfaces with mean curvature $\vec{h}$ and normal velocity $\vec{\tilde{v}}$, and 
$\vec{v}$ satisfy \eqref{Brakke ineq} with $V_t(\phi)=\int_{M_t} \phi (x,T_x M_t) \, {\rm d} \mathcal{H}^{d-1} (x)$ for $\phi \in C_c (\mathbf{G}_{d-1}(\R^d ))$, where $T_{x}M_t$ is the tangent hyperplane to $M_t$ at $x$ and $\mathbf{G}_{d-1}(\R^d)$ is the Grassmannian bundle $\R^d\times \mathbf{G}(d,d-1)$, where $\mathbf{G}(d,d-1)$ is the Grassmannian of unoriented $(d-1)$-planes. Then $\vec{\tilde{v}}$ satisfies \eqref{Brakke ineq} with equality
    \begin{align}
        \|V_t\|(\phi(\cdot,t))\Big|_{t=t_1}^{t_2}&=\int_{t_1}^{t_2}\int_{\R^d} \left(\phi_t+(\nabla\phi-\phi\vec{h})\cdot\vec{\widetilde{v}}\right)\,{\rm d}\|V_t\|\dt.\label{eq brakke inequality classical v}
    \end{align}
    Hence, subtracting \eqref{eq brakke inequality classical v} to \eqref{Brakke ineq} and differentiating in time,  we get
    \begin{equation*}
        0\leq \int_{\R^d} \left((\nabla\phi-\phi\vec{h})\cdot(\vec{v}-\vec{\widetilde{v}})\right)\,{\rm d}\|V_t\|
    \end{equation*}
    for every $\phi\in C^1_c({\R^d};[0,\infty))$ and every $t\in (0,T)$. In particular, fix $\phi\in C^1_c({\R^d};[0,\infty))$, $x_0 \in M_t$, and let $\phi_r(x)\coloneqq \frac{1}{r^{d-2}}\phi\left(\frac{x{-}x_0}{r}\right)$ for $r>0$: The last inequality for the test function $\phi_r$ becomes, by change of variables and sending $r\rightarrow 0^+$, 
    \begin{equation*}
        0\leq \int_{T_{x_0}M_t} \left(\nabla\phi{\cdot}(\vec{v}-\vec{\widetilde{v}})\right)\dH{d-1}.
    \end{equation*}
    Since $x_0\in M_t$ and $\phi\in C^1_c({\R^d}\times(0,T);[0,\infty))$ are arbitrary, then $\vec{v}=\vec{\widetilde{v}}$, and so \eqref{Brakke ineq} implies that $\vec{v}$ coincides with the classical normal velocity in the smooth setting.

In the case of the volume preserving mean curvature flow equation \eqref{vpmcf}, however, Brakke inequality in the classical form \eqref{Brakke ineq} poses significant challenges to the analysis. Indeed, for $\vec{v}=\vec{h}-\lambda\vec{\nu}$,
\eqref{Brakke ineq} can be rewritten as
    \begin{align}
       \|V_t\|(\phi(\cdot,t))\Big|_{t=t_1}^{t_2}&\leq\int_{t_1}^{t_2}\int_{\R^d} \left(\phi_t+\nabla\phi{\cdot}\left(\hh-\lambda \vec{\nu}\right)-\phi\left(|\hh|^2-\lambda\hh{\cdot}\vec{\nu}\right)\right)\,{\rm d}\|V_t\|\dt\notag\\
        &=\int_{t_1}^{t_2}\int_{\R^d} \left(\phi_t+\nabla\phi{\cdot}\left(\hh-\lambda \vec{\nu}\right)-\phi\left(\frac{|\hh|^2}{2}+\frac{|\hh-\lambda\vec{\nu}|^2}{2}-\frac{|\lambda|^2}{2}\right)\right)\,{\rm d}\|V_t\|\dt\notag\\
        &=\int_{t_1}^{t_2}\int_{\R^d}\left(\phi_t+\nabla\phi{\cdot}\left(\hh-\lambda \vec{\nu}\right)-\phi\left(\frac{|\hh|^2}{2}+\frac{|\vec{v}|^2}{2}\right)+\phi\frac{|\lambda|^2}{2}\right)\,{\rm d}\|V_t\|\dt,\label{eq brakke inequality correspondence}
    \end{align}
which is analogous to the monotonicity formula obtained in \cite{LM}.    
Such inequality is usually obtained by approximation, with some sequences $h^\eps$, $v^\eps$, and $\lambda^\eps$  weakly converging (in a suitable sense) to $\vec{h}$, $\vec{v}$, and $\lambda$, respectively. The terms involving $|\vec{h}|^2$ and $|\vec{v}|^2$ can be hence obtained by the weak lower-semicontinuity of the sequence of measure-function pairs (cf. \cite{Hutchinson}). The $|\lambda|^2$ term, however, has the wrong sign and is not amenable to such analysis. In particular, the estimate 
\begin{equation*}
\limsup_{\eps\rightarrow0}\int_{t_1}^{t_2}\int_{\R^d}\phi |\lambda^\eps|^2\,{\rm d}\|V_t^\eps\|\dt\leq \int_{t_1}^{t_2}\int_{\R^d}\phi |\lambda|^2\,{\rm d}\|V_t\|\dt,
\end{equation*}
for a suitable approximating varifold $\{V^\eps_t\}_{t\in[0,T]}$, holds only assuming $\lambda^\eps\rightarrow\lambda$ in the strong $L^2(0,T)$-topology, which cannot be expected in general. 

In this paper, we thus show the existence of a varifold flow satisfying the following inequality for some $C>0$ (cf. Definition \ref{def Brakke flow}), 
\begin{align*}
   \|V_t\|(\phi(\cdot,t))\Big|_{t=t_1}^{t_2}\leq&\int_{t_1}^{t_2}\int_{\R^d} \left(\phi_t+\nabla\phi{\cdot}\left(\hh-\lambda \vec{\nu}\right)-\phi\left(\frac{|\hh|^2}{2}+\frac{|\vec{v}|^2}{2}\right)\right)\,{\rm d}\|V_t\|\dt\notag\\
    &+r^{d-1}C(1+t_2-t_1)\|\phi\|_{L^\infty}
\end{align*}
for every $r\in (0,1)$, for every $0\leq t_1< t_2\leq T$, and for every $\phi\in C^1_c({\R^d}\times[0,T];[0,\infty))$ with $\operatorname{spt}{\phi}\subset \subset B_r(x_0)$, for some $x_0\in {\R^d}$. 
This inequality coincide with the classical Brakke inequality \eqref{eq brakke inequality correspondence} except for the term $r^{d-1}C(1+t_2-t_1)\|\phi\|_{L^\infty}$, which corresponds to the $|\lambda|^2$-term in the last line of \eqref{eq brakke inequality correspondence}. More precisely, we can prove
    \begin{equation*}
        r^{d-1}C(1+t_2-t_1)\|\phi\|_{L^\infty}\geq \limsup_{\eps\rightarrow0}\int_{t_1}^{t_2}\int_{\R^d}\phi |\lambda^\eps|^2\,{\rm d} \|V_t^\eps\|\dt,
    \end{equation*}
    cf. \eqref{bound L^2 norm lambda^eps } and the following part of the proof of Theorem \ref{thm: existence Brakke}.
Moreover, this modified inequality still characterises the normal velocity in the smooth setting, since, reasoning as above with the test function $\phi_r(x)\coloneqq \frac{1}{r^{d-2}}\phi\left(\frac{x{-}x_0}{r}\right)$, the additional term $r^{d-1}C(1+t_2-t_1)\|\phi\|_{L^\infty}={\rm O}(r)$ as $r\rightarrow 0$.

The paper is structured as follows. In Section \ref{sec: def and hp}, we recall the definition of $L^2$-flow, define the Brakke notion of solution we are interested in, and provide the assumption on the initial datum $U_0$. In Section \ref{sec: mai results}, we state our main results: In Section \ref{sec: mai results L^2}, we recollect the existence result for the $L^2$-flow of \cite{TakasaoVPMCF}, whereas, in Section \ref{sec: mai results Brakke}, we present and prove the existence of a Brakke-flow.

\section{Definitions and assumptions}\label{sec: def and hp}

\subsection{Notions of solution: \texorpdfstring{$L^2$}{}-flow and Brakke-flow}
Let $\Omega$ be an open subset of $\R^d$ or the $d$-dimensional torus $\mathbb{T}^{d} =(\R /\Z)^d$. In this section we recall the definition of solution to \eqref{vpmcf} in the sense of $L^2$-flow and provide the Brakke notion of solution we are interested in. Let $\{\mu_t\}_{t\in[0,T)}$ be a family of Radon measures on $\Omega$.

\begin{definition}[$L^2$-flow {\cite[Def. 2.1]{TakasaoVPMCF}}]\label{def L^2 flow}
    We say that $\{\mu_t\}_{t\in[0,T)}$ is an {\rm $L^2$-flow} with generalized velocity vector $\vv$ if:
    \begin{enumerate}
        \item for a.e. $t\in (0,T)$, there exists a $(d{-}1)$-integral varifold $V_t\in \V_{d-1}(\Omega)$ such that $\mu_t=\|V_t\|$ and having generalized mean curvature vector $\hh\in L^2(\m;\R^d)$;
        \item $\vv\in L^2(0,T;L^2(\m;\R^d))$ and 
        \begin{equation}\label{def L^2 flow v perpendicular}
            \vv(x,t)\perp T_x\m \qquad \text{ for $\mu$-a.e. } (x,t)\in \Omega\times (0,T),
        \end{equation}
        where $\dmu\coloneqq \dmu_t\dt$ and $T_x\m\in \mathbf{G}(d,d{-}1)$ is the approximate tangent space of $V_t $ at $x$;
        \item there exists a constant $C_T>0$ such that 
        \begin{equation}\label{def L^2 flow v weak velocity}
            \left|\int_0^T\int_\Omega (\eta_t+\nabla\eta{\cdot}\vv)\dmu_t\dt\right|\leq C_T\|\eta\|_{L^\infty(\Omega\times(0,T))}
        \end{equation}
        for every $\eta\in C^1_c(\Omega {\times} (0,T))$.
    \end{enumerate}
\end{definition}

\begin{remark}
    Notice that if $\vv$ is the velocity of a smooth hypersurface $M_t\subset U$, then 
    \begin{equation*}
        \frac{\rm d}{{\rm d}t}\int_{M_t}\eta \,  \dhx{d-1}=\int_{M_t}(\nabla^\perp \eta -\eta \hh){\cdot}\vv+\eta_t \, \dhx{d-1},
    \end{equation*}
    for every $\eta \in C^1_c(\Omega{\times}(0,T))$, where $\nabla^\perp\eta\coloneqq (\nabla \eta{\cdot}\nu)\nu$. If $\int_0^T\int_{M_t}|\vv{\cdot}\hh|\dhx{d-1}\dt<\infty$ and $\vv$ satisfies \eqref{def L^2 flow v perpendicular} (i.e., is perpendicular to $T_xM_t$), then \eqref{def L^2 flow v weak velocity} follows. Conversely, one can show that if $\vv$ satisfies \eqref{def L^2 flow v perpendicular} and \eqref{def L^2 flow v weak velocity}, then it is the normal velocity vector of $M_t$. In this sense $\vv$ 
    in Definition \ref{def L^2 flow} is the weak normal velocity vector of $\mu_t$.
\end{remark}

For a set $A\subset \Omega$ of finite perimeter, we denote the reduced boundary of $A$ by $\partial ^\ast A$.
\begin{definition}[Volume-preserving Brakke-flow]\label{def Brakke flow}
  Let $T \in (0,\infty)$, $\Omega =\R^d$ or $\mathbb{T}^d$, and $U_0 \subset \Omega$ be a set of finite perimeter. We say that $(\{\mu_t\}_{t\in[0,T]}, \psi, \lambda)$ is a {\rm volume-preserving Brakke-flow} with initial data $\partial ^\ast U_0$ if
  $\mu_t$ is a Radon measure on $\Omega$ for any $t \in [0,T]$, $\psi\in BV(\Omega \times [0,T];\{0,1\})\cap C ^{\frac12} ([0,T]; L^1 (\Omega))$, $\lambda \in L^2 (0,T)$, and they are such that:
\begin{enumerate}
\item it holds that $\mu_0 = \mathcal{H}^{d-1} \lfloor_{\partial ^\ast U_0}$ and there exists $C>0$ such that $\sup_{r>0, x_0 \in \Omega} \frac{\mu_0 (\overline{B_r (x_0)})}{r^{d-1}}<C$;
\item If $\Omega=\R^d$, then $\{\mu_t\}_{t \in [0,T]}$ satisfies
$\sup_{t \in [0,T]} \mu_t (K) <\infty$ for any compact set $K \subset \subset \Omega$. If $\Omega=\mathbb{T}^d$, then $\{\mu_t\}_{t \in [0,T]}$ satisfies
$\sup_{t \in [0,T]} \mu_t (\mathbb{T}^d) <\infty$;
\item for a.e. $t\in [0,T]$, there exists a $(d{-}1)$-integral varifold $V_t\in \V_{d-1}(\Omega)$ such that $\mu_t=\|V_t\|$ and having generalized mean curvature vector $\vec{h} \in L^2((0,T); L^2 _{\mathrm{loc}}(\mu_t;\R^d))$;
\item 
$\psi$ satisfies the volume preserving property, i.e.,
\[
\int _{\Omega} \psi (x,t) \, \dx  = \mathcal{L}^{d} (U_0)
\qquad \text{for all} \ t \in [0,T);
\]
\item for every $t \in [0,T]$ and for any nonnegative $\phi \in C_c (\Omega\times [0,T] ;[0,\infty))$, 
we have 
\begin{equation}\label{psi absolutely cont wrt mu}
 \| \nabla \psi (\cdot ,t) \| (\phi (\cdot,t) ) \leq \mu _t (\phi (\cdot,t));
\end{equation}
 \item there exists a constant $C$ such that for any $r \in (0,1)$ it holds that
\begin{align}
    \mu_t(\phi(\cdot,t))\Big|_{t=t_1}^{t_2}\leq&\int_{t_1}^{t_2}\int_{\Omega} \left(\phi_t+\nabla\phi{\cdot}\vec{v}-\phi\left(\frac{|\hh|^2}{2}+\frac{|\vec{v}|^2}{2}\right)\right)\dmu_t\dt\notag\\
    &+r^{d-1}C(1+t_2-t_1)\|\phi\|_{L^\infty}\label{Brakke ineq weak}
\end{align}
for every $0\leq t_1< t_2\leq T$ and for every $\phi\in C^1 _c({\R^d}\times[0,T];[0,\infty))$ with $\operatorname{spt}{\phi}(\cdot, t)\subset \subset B_r(x_0)$ for every $t\in [0,T]$, for some $x_0 \in \R^d$
(when $\Omega =\mathbb{T}^d$, $\phi$ is periodically extended to be included in $C ^1 (\mathbb{T}^d\times [0,T];[0,\infty))$). Here, $\vec{v} \in L^2 ((0,T); L^2 _{\mathrm{loc}} (\mu_t ;\R^d))$ satisfies $\vec{v}:= \vec{h} 
- \lambda \frac{{\rm d} \| \nabla \psi (\cdot,t) \|}{\dmu _t} \vec{\nu}$ in the sense of distributions, namely, 
 for every $\vec{\phi}\in C_c({\Omega}\times(0,T);\R^d)$, it holds that 
 \begin{equation}\label{v=h-lambda distrib}
     \int_{0} ^T \int _{\Omega} \vec{v} \cdot \vec{\phi} \, {\rm d} \mu_t {\rm d} t
=\int_{0} ^T \int _{\Omega} \left( \vec{h} 
- \lambda \frac{{\rm d} \| \nabla \psi (\cdot,t) \|}{\dmu _t} \vec{\nu} \right) \cdot \vec{\phi} \,  \dmu_t \dt,
 \end{equation}
 where $\frac{d \| \nabla \psi (\cdot,t) \|}{d\mu _t}$ is the Radon-Nikodym derivative and $\vec \nu$ is the inner unit normal vector of $\{ \psi (\cdot,t) =1 \}$
 on $\operatorname{spt} \| \nabla \psi (\cdot,t) \|$.
\end{enumerate}
\end{definition}


We notice that, analogously to \cite[Thm. 4.3]{StuvardTonegawa}, any volume-preserving Brakke flow is an $L^2$-flow with generalized velocity vector  ``$\vec{v}=\vec{h}-\lambda \frac{{\rm d}\|\nabla \psi(\cdot,t)\|}{\dmu_t} \vec{\nu}\,$". Namely, under suitable assumptions, which can be considered without loss of generality (cf. \ref{(a)} in Theorem \ref{mainthm1}), we have the following result.

\begin{proposition}
    Let $T \in (0,\infty)$, $\ \Omega=\R^d$ or $\mathbb{T}^d$ and  $(\{\mu_t\}_{t\in[0,T)}, \psi,\lambda)$ be a volume-preserving Brakke-flow with $\mu_t(\Omega)\leq \mu_0(\Omega)<\infty$
    for almost every $t\in [0,T]$. Suppose that $U\subset \Omega$ be a bounded open set. Then $\{\mu_t\}_{t\in[0,T)}$ satisfies 
    \begin{equation}\label{Lemma brakke flow v = h}
            \left|\int_0^T\int_U \left(\eta_t+\nabla\eta{\cdot}\!\left(\hh-\lambda \frac{{\rm d}\|\nabla \psi(\cdot,t)\|}{\dmu_t}\vec{\nu}\right)\right)\!\dmu_t\dt\right|\leq C_T\|\eta\|_{L^\infty(U\times(0,T))}
        \end{equation}
    for every $\eta\in C^1_c({U}{\times} (0,T))$, for some constant $C_T>0$. In particular, $\{\mu_t\}_{t\in[0,T)}$ is a $L^2$-flow with generalized velocity vector $\vv=\hh-\lambda \frac{{\rm d}\|\nabla \psi(\cdot,t)\|}{\dmu_t}\vec{\nu}$ on $U$. 
\end{proposition}

\begin{proof}
We show only the case $\Omega =\R^d$ (the case $\Omega =\mathbb{T}^d$ is analogous).
 Notice that, for a.e. $t \in (0,T)$, $V_t$ is a $(d-1)$-integral varifold. For such $t$, 
 \begin{equation}\label{perp1}
 \vec{h} (x,t) \perp T_x \mu_t \quad \text{for} \quad \mu_t\text{-a.e.} \ x 
 \end{equation}
 by the perpendicularity of mean curvature \cite[Chap. 5]{Brakke}.
 Choose $t \in (0,T)$ such that \eqref{perp1} and $|\lambda (t)|<\infty$ hold.
 Let $M_t$ and $\tilde M_t$ be $\mathcal{H}^{d-1}$-measurable countably $(d-1)$-rectifiable sets such that $\mu_t = \theta_t \mathcal{H}^{d-1} \lfloor_{M_t}$ and 
 $\| \nabla \psi (\cdot,t) \| = \mathcal{H}^{d-1} \lfloor_{\tilde M_t}$, where $\theta_t$ is the multiplicity function of $\mu_t$. To show 
 \eqref{def L^2 flow v perpendicular}, we only need to prove it for $\mathcal{H}^{d-1}$-a.e. $x\in M_t$. Note that we may choose $\tilde M_t =\partial^\ast \{ x \mid \psi (x,t)=0 \}$, and \eqref{perp1} holds for $\mathcal{H}^{d-1}$-a.e. $x \in M_t$, since if $\mu_t (N)=0$ then $\mathcal{H}^{d-1} (M_t \cap N)=
 \lim_{j\to \infty}\mathcal{H}^{d-1} (M_t^j \cap N)\leq \lim_{j\to \infty} j \mu_t (N)=0$, where $M_t ^j= \{ x \in M_t \mid \theta _t (x)>1/j \}$.
 For $\mathcal{H}^{d-1}$-a.e. $x \in M_t \cap \tilde{M}_t$, $T_x \mu _t =T_x M_t = T_x \tilde M_t$ and hence $\vec\nu (x,t) \perp T_x \mu_t$ for such $x$ (see \cite[Proposition 2.85]{Ambrosio}). By this and \eqref{perp1}, \eqref{def L^2 flow v perpendicular} holds for $\mathcal{H}^{d-1}$-a.e. $x \in M_t \cap \tilde{M}_t$.
 In addition, for $\mathcal{H}^{d-1}$-a.e. $x \in M_t \setminus \tilde M_t$, 
 $\lim_{r\downarrow 0}\frac{\mu_t (B_r (x))}{r^{d-1}} >0$ and $\lim_{r\downarrow 0}\frac{\mathcal{H}^{d-1} (\tilde M_t \cap B_r (x))}{r^{d-1}} =0$ (see \cite[Theorem 3.5]{Simon}), and hence
 \[
 \frac{{\rm d}\| \nabla \psi (\cdot,t) \|}{{\rm d} \mu _t }(x)
 =\lim_{r\downarrow 0} 
 \frac{\mathcal{H}^{d-1} (\tilde M_t \cap B_r (x))/r^{d-1}}{\mu _t (B_r (x))/r^{d-1}}
 =0.
 \] 
 Therefore the term $\lambda \frac{{\rm d}\|\nabla \psi(\cdot,t)\|}{\dmu_t}\vec{\nu}$ vanishes for
 $\mathcal{H}^{d-1}$-a.e. $x \in M_t \setminus \tilde M_t$ and thus \eqref{def L^2 flow v perpendicular} holds for
 $\mathcal{H}^{d-1}$-a.e. $x \in M_t$.
 
 Hence, it remains to show \eqref{Lemma brakke flow v = h}. 
 Assume that $\phi\in C^2 _c({\R^d}\times (0,T);[0,\infty))$ with $\operatorname{spt}{\phi}(\cdot, t)\subset \subset B_r(x_0)$ for every $t\in [0,T]$, for some $x_0\in {\R^d}$. Then by \eqref{Brakke ineq weak}, we have
 \begin{equation*}
 \begin{split}
 & \int_0 ^T \int_{\R^d} \phi \frac{|\vec v|^2}{2} \, \dmu_t \dt \\
 \leq & \, \| \phi_t \|_\infty \sup_{t \in [0,T]} \mu_t (\mathrm{spt} \, \phi (\cdot,t) ) T
 + \int_0 ^T \int_{\R^d} 
 \frac{|\nabla \phi|^2}{\phi} + \frac{\phi |\vec v|^2}{4} \, \dmu_t \dt
 + r^{d-1} C(1+T) \| \phi \|_\infty,
 \end{split}
 \end{equation*}
 where we used Young's inequality. 
 Note that $\sup_{t \in [0,T]} \mu_t (\mathrm{spt} \, \phi (\cdot,t)) \leq \sup_{t \in [0,T]} \mu_t (\overline{B_r (x_0)}) <\infty$ by (2) of Definition \ref{def Brakke flow}.
 Since there exists $C>0$ such that 
 $\| \frac{|\nabla \zeta|^2}{\zeta} \|_\infty \leq C \| \zeta  \|_{C^2}$ for any $\zeta \in C_c ^2 (\R^d)$ \cite[Lemma 3.1]{TonegawaBook}, we obtain
 \begin{equation}\label{L2-estimate-velo}
 \begin{split}
 & \int_0 ^T \int_{\R^d} \phi |\vec v|^2 \, \dmu_t \dt \\
 \leq & \, C \left( (\| \phi_t \|_\infty+\|\phi \|_{C^2}) \sup_{t \in [0,T]} \mu_t (\overline{B_r (x_0)} ) T
 + r^{d-1} (1+T) \| \phi \|_\infty\right).
 \end{split}
 \end{equation} 
 For every $r>0$, $x_0\in \R^d$, and every nonnegative test function $\phi\in C^1_c({\R^d}{\times} (0,T);[0,\infty))$ with $\operatorname{spt}\phi(\cdot,t)\subset \subset B_r(x_0)$ for every $t\in [0,T]$, by \eqref{Brakke ineq weak} we have 
\begin{equation}\label{eq brakke lemma brakke imply L^2}
\begin{split}
     0\leq & \, \int_{0}^{T}\int_{\R^d} \left(\phi_t+\nabla\phi{\cdot}\left(\hh-\lambda \frac{{\rm d}\|\nabla \psi(\cdot,t)\|}{\dmu_t}\vec{\nu}\right)\right)\dmu_t\dt+r^{d-1}C(1+T)\|\phi\|_{L^\infty}
     +\mu_0 (\phi (\cdot,0))\\
     \leq & \,
     \int_{0}^{T}\int_{\R^d} \left(\phi_t+\nabla\phi{\cdot}\left(\hh-\lambda \frac{{\rm d}\|\nabla \psi(\cdot,t)\|}{\dmu_t}\vec{\nu}\right)\right)\dmu_t\dt+r^{d-1}C(1+T)\|\phi\|_{L^\infty},
    \end{split}
 \end{equation}
where we used the upper bound of the density of the initial data ((2) of Definition \ref{def Brakke flow}), namely, 
$\mu_0 (\mathrm{spt} \, \phi (\cdot,0)) \leq \mu_0 (\overline{B_r (x_0)}) \leq r^{d-1} C$. 
Then, let $\eta\in C^1_c({\R^d}{\times} (0,T))$ be a non-zero test function with $\operatorname{spt}\eta\subseteq B_{r/2}(x_0)\times [\delta,T]$ for some $\delta>0$ and some fixed $r>0$. Moreover, let $\zeta_{\delta,r}\coloneqq \zeta_\delta\zeta_r \in C^3 _c(B_r(x_0)\times (0,T])$  be cut-off functions such that 
\begin{enumerate}[label=\roman*)]
    \item $0\leq \zeta_\delta,\zeta_r\leq 1$;
    \item $\zeta_\delta\equiv 1$ on $[\delta,T]$ and $\zeta_\delta\equiv 0$ on $(0,\delta/2]$;
    \item $\zeta_r\equiv 1$ on $B_{r/2}(x_0)$ and $\zeta_r\equiv 0$ on $\R^d\setminus B_{r}(x_0)$;
    \item $\| \zeta_\delta'\|_{L^\infty}\leq  c/\delta$ and $\| \nabla\zeta_r\|_{L^\infty}\leq  c/r$ for some $c>0$.
\end{enumerate}
Thus, we have 
\begin{equation*}
    -\zeta_{\delta,r}\leq \frac{\eta}{\|\eta\|_{L^\infty}}\leq \zeta_{\delta,r}
\end{equation*}
and \eqref{eq brakke lemma brakke imply L^2} applied to $\phi\coloneqq\zeta_{\delta,r}\|\eta\|_{L^\infty}\pm\eta\geq 0$ gives
\begin{align*}
    &\left|\int_{0}^{T}\int_{\R^d} \left(\eta_t+\nabla\eta{\cdot}\left(\hh-\lambda \frac{{\rm d}\|\nabla \psi(\cdot,t)\|}{\dmu_t}\vec{\nu}\right)\right)\dmu_t\dt\right|
    \\
    &\quad\leq \left(\int_{0}^{T}\int_{\R^d} \left((\zeta_{\delta,r})_t+\nabla\zeta_{\delta,r}{\cdot}\left(\hh-\lambda \frac{{\rm d}\|\nabla \psi(\cdot,t)\|}{\dmu_t}\vec{\nu}\right)\right)\dmu_t\dt+2C(1+T)\right)\|\eta\|_{L^\infty}\\
&\quad\stackrel{\eqref{v=h-lambda distrib}}{\leq } \left(\int_{0}^{T}\int_{\R^d} \Big((\zeta_{\delta,r})_t+\nabla\zeta_{\delta,r}{\cdot}\vv\Big)\dmu_t\dt+2C(1+T)\right)\|\eta\|_{L^\infty}.
\end{align*}
In order to bound the right-hand side, notice that, by the assumptions on $\zeta_{\delta,r}$, 
\begin{align*}
    & \int_{0}^{T}\int_{\R^d} \Big((\zeta_{\delta,r})_t+\nabla\zeta_{\delta,r}{\cdot}\vv\Big)\dmu_t\dt \\
    \leq &\, \frac{c}{\delta}\int_{\delta/2}^\delta \mu_t(B_r(x_0))\dt+ \int_0^T\int_{\R^d} |\nabla \zeta_r| \zeta_\delta |\vv|\dmu_t\dt\\
    \leq & \, c \sup_{t \in [0,T]}\mu_t (\overline{B_r (x_0)})+\sqrt{T}\sqrt{\sup_{t \in [0,T]}\mu_t (\overline{B_r (x_0)})}\left(\int_0^T\int_{\R^d} |\nabla \zeta_r|^2 \zeta_\delta ^2|\vv|^2\dmu_t\dt\right)^{1/2}.
\end{align*}
By \eqref{L2-estimate-velo}
and $|\nabla \zeta_r|^2 \zeta_\delta^2 \in C_c ^2 (\R^d \times (0,T);\R^d_+)$, the right hand side is bounded independently of $\delta$ (the bound depends on $r$, which we have however fixed above).
Consequently, for every $r\in (0,1)$ and $x_0\in \R^d$, \eqref{Lemma brakke flow v = h} holds for every $\eta\in C^1_c(U\times (0,T))$ with $\operatorname{spt}\eta(\cdot,t)\subseteq B_{r/2}(x_0)$ and  by partitions of unity and the boundedness of $U$ we conclude. 
\end{proof}


\subsection{Assumptions on the initial data}
Henceforth, let $\Omega=\mathbb{T}^d=(\R/\Z)^d$. In this paper we consider the periodic boundary condition for \eqref{vpmcf}.
Let $U _0 \subset \subset \Omega$ be a bounded open set with the following properties.
\begin{enumerate}
\item There exists $D_0 >0$ such that
\begin{equation}\label{U01}
\sup _{x \in \Omega , 0< r <1} \frac{\HH^{d-1} (M_0 \cap B_r (x))}{\omega _{d-1} r^{d-1}}
\leq D_0,
\end{equation}
where $M_0 = \partial U_0$.
\item There exists a family of open sets $\{ U_0 ^i\}_{i=1} ^\infty$ such that
$U_0 ^i$ has a $C^3$ boundary 
$M_0 ^i = \partial U_0 ^i$ for any $i $ and the following hold:
\begin{equation}\label{U02}
\lim_{i\to \infty}
\mathcal{L}^d (U_0 \triangle U_0 ^i)=0
\qquad
\text{and}
\qquad
\lim_{i\to \infty}
\| \nabla \mathbbm{1} _{U_0 ^i} \|= \| \nabla \mathbbm{1} _{U_0} \|
\quad \text{as Radon measures}.
\end{equation}
\end{enumerate}

Moreover, we define a periodic function $\varphi _0 ^{\eps} \in C^3 (\Omega)$ as a suitable smoothing of $2\mathbbm{1} _{U_0}{-}1$, see \cite{TakasaoVPMCF} for more details.

\section{Main results}\label{sec: mai results}
In this section, we assume $\Omega:=\mathbb{T}^d$.
Let us consider the following Allen-Cahn equation as in \cite{ Golovaty,TakasaoVPMCF}
\begin{equation}
\left\{ 
\begin{array}{ll}
\eps \varphi ^{\eps} _t =\eps \Delta \varphi ^{\eps} -\dfrac{W' (\varphi ^{\eps})}{\eps }+ \lambda ^{\eps} \sqrt{2W(\varphi ^\eps)} ,& (x,t)\in \Omega \times (0, \infty),  \\
\varphi ^{\eps} (x,0) = \varphi _0 ^{\eps} (x) ,  &x\in \Omega,
\end{array} \right.
\label{ac}
\end{equation}
where the non-local term $\lambda^{\eps}$ is given by
\begin{equation}\label{lambdadef}
\lambda^{\eps}(t) =\frac{1}{\eps ^\alpha} 
\left( \int_{\Omega} k(\varphi ^\eps _0 (x))\, \dx 
- \int_{\Omega} k(\varphi ^\eps (x,t))\, \dx \right)
\end{equation}
with $\alpha\in(0,1)$, $W$ is the standard double-well potential given by $W(s)\coloneqq (1-s^2)^2/2$, and  $k(s)\coloneqq \int_0^s\sqrt{2W(r)}\dr=s-\frac{1}{3}s^3$. In this section, we assume that the initial data $\varphi_0 ^\eps$ satisfies $\sup_{x \in \Omega} |\varphi _0 ^\varepsilon (x)| <1$ for any $\varepsilon \in (0,1)$, so that $\sup_{(x,t) \in \Omega \times (0,T)} |\varphi ^\varepsilon (x,t)| <1$ for any $\varepsilon \in (0,1)$ and $T>0$, by the maximum principle.
We remark that \cite{BMR} recently  proposed a phase field model similar to \eqref{ac} and studied the asymptotic behavior of its solution.
Note that under suitable assumptions on $\varphi _0 ^\eps$,
 global existence and uniqueness of 
the solution to \eqref{ac} follow by the standard PDE theory.
Equation \eqref{ac} is the gradient flow of the energy $E^\eps$, given by 
\begin{align*}
E^{\eps} (t) \,
& \coloneqq   E_S ^\eps (t) + E_P ^\eps (t),
\end{align*}
where
\begin{align*}
E_S ^\eps (t)&\coloneqq
\,\int _{\Omega} 
\left( \frac{\eps \vert \nabla \varphi ^\eps (x,t) \vert ^2}{2} 
+ \frac{W(\varphi ^\eps (x,t))}{\eps} \right) \, \dx,\\
E_P ^\eps (t)&\coloneqq
\frac{1}{2\eps ^\alpha}
\left( \int_{\Omega} k(\varphi ^\eps _0 (x))\, \dx 
- \int_{\Omega} k(\varphi ^\eps (x,t))\, \dx \right)^2.
\end{align*}
The singular limit of its gradient flow \eqref{ac} is the volume preserving mean curvature flow in the sense of $L^2$-flow \cite{TakasaoVPMCF}, where $E^\eps_S$ corresponds to the surface energy, whereas $E^\eps_P$ penalises the difference of the volumes.
Let also $\mu_t^\eps$ be the Radon measure defined as 
\begin{equation}\label{mu^eps}
    \mu_t^\eps(\phi)\coloneqq \frac{1}{\sigma}\int_\Omega \phi\left(\frac{\eps |\nabla \varphi^\eps(x,t)|^2}{2}+\frac{W(\varphi^\eps(x,t))}{\eps}\right)\dx, \qquad \forall \phi\in C_c(\Omega)
\end{equation}
 for $t\in[0,T)$, where $\sigma\coloneqq \int_{-1}^1\sqrt{2W(s)}\ds$, and $\mu^\eps$ be defined by ${\rm d} \mu^\eps= {\rm d} \mu^\eps_t\dt$. 
We denote the approximate velocity vector $\vec{v} ^{\, \varepsilon }$ and the approximate mean curvature $h^\eps$ by
\begin{equation}\label{approximate v and h}
\vec{v} ^{\, \varepsilon } \coloneqq \left\{ 
\begin{array}{ll}
\dfrac{- \varphi ^{\varepsilon } _t }{\vert \nabla \varphi ^{\varepsilon} \vert }
\dfrac{\nabla \varphi ^{\varepsilon }}{\vert \nabla \varphi ^{\varepsilon } \vert }, & \text{if} \ 
\vert \nabla \varphi ^{\varepsilon } \vert \not=0, \\
\qquad 0 , & \text{otherwise},
\end{array} \right.\qquad h^\eps\coloneqq \Delta \varphi^\eps-  \frac{W'(\varphi^\eps)}{\eps^2},
\end{equation}
respectively.
In the following sections we often consider subsequences $(\eps_i)_{i\in \N}$ but write $\eps$ for notational simplicity. Let us first recall that, for $\varphi ^{\eps} _0$ and $\mu _0 ^{\eps}$, we have the following properties.
\begin{proposition}[{\cite[Prop. 2]{TakasaoVPMCF}}]\label{prop2.2}
Let $\{\varepsilon _i \}_{i=1} ^\infty$ be a positive sequence with $\varepsilon_i \to 0$ as $i\to \infty$.
For suitable $\varphi^{\eps_i} _0$ given as in \cite{TakasaoVPMCF}, up to a nonrelabeled subsequence the following hold:
\begin{enumerate}
\item For any $i\geq 1$ and $x \in \Omega$, we have
$\displaystyle
\frac{\eps_i \vert \nabla \varphi^{\eps_i} _0 (x) \vert ^2 }{2} 
\leq \frac{W(\varphi ^{\eps_i} _0 (x))}{\eps_i} 
$.
\item
There exists $D_1= D_1 (D_0) >0$ such that
\begin{equation}
\max \left\{
\sup _{i\geq 1} \mu _0 ^{\eps_i} (\Omega) 
, \sup _{i\geq 1, \ x \in \Omega, \ r \in (0,1)} 
\frac{\mu _0 ^{\eps_i} (B_r (x))}{\omega ^{d-1} r^{d-1}}
\right\}
\leq D_1.
\label{d1}
\end{equation}
\item $\mu _0 ^{\eps_i} \to \HH^{d-1} \lfloor _{M_0}$ as Radon measures, that is,
\[
\int _{\Omega} \phi \, {\rm d} \mu _0 ^{\eps_i} \to \int _{M_0} \phi \, {\rm d}\HH^{d-1}
\qquad 
\text{for any} \ \phi \in  C_c(\Omega).
\]
\item For $\psi ^{\eps_i}_0 =\frac12 ( \varphi ^{\eps_i}_0 +1)$,
$\lim _{i\to \infty} \psi ^{\eps_i}_0 = \mathbbm{1} _{U_0} $ in $L^1$
and $\lim _{i\to \infty} \| \nabla \psi_0 ^{\eps_i } \| = \| \nabla \mathbbm{1} _{U_0} \| $
as Radon measures.
\end{enumerate}
\end{proposition}

\subsection{Previous results: Existence of an \texorpdfstring{$L^2$}{}-flow}\label{sec: mai results L^2}
The existence of an $L^2$-flow solving \eqref{vpmcf} has been shown in \cite{TakasaoVPMCF}. In particular we have the following results.

To simplify notation, the subscript $i$ in $\varepsilon_i$ is sometimes omitted. In such cases, the expression $\varepsilon \to 0$ means taking the limit as $i \to \infty$ for a positive sequence $\{\varepsilon_i\}_{i=1} ^\infty$ such that $\varepsilon_i \to 0$.

\begin{theorem}[{\cite[Thm. 3, Thm. 14, Lemma 1]{TakasaoVPMCF}}]\label{mainthm1}
Suppose that $d\geq 2$ and $U_0$ satisfies \eqref{U01} and \eqref{U02}.
For any $\eps>0$, let $\varphi _0 ^{\eps}$ be defined so that all the claims of Proposition \ref{prop2.2} are satisfied
and
$\varphi ^{\eps} $ be a solution to
\eqref{ac} with initial data $\varphi _0 ^{\eps}$.
Then, up to a nonrelabeled subsequence, the following hold:
\begin{enumerate}[label=\textnormal{(\alph*)}]
\item\label{(a)} There exist a countable subset $B \subset [0,T)$ and
a family of Radon measures $\{\mu _t\}_{t \in [0,T)}$
on $\Omega$ such that
\[
\mu _0 = \HH^{d-1} \lfloor_{M_0}, 
\qquad 
\mu _t ^{\eps} \to \mu _t \ \ \ \text{as Radon measures for any} \ t \in [0,T),
\]
\[
\mu _s (\Omega) \leq \mu _t (\Omega) \qquad \text{for any} \ s,t \in [0,T) \setminus B 
\ \text{with} \  0 \leq t <s <T,
\]
and for a.e. $t\in [0,T]$, there exists a $(d{-}1)$-integral varifold $V_t\in \V_{d-1}(\Omega)$ such that $\mu_t=\|V_t\|$ and having generalized mean curvature vector $\vec{h} (\cdot,t) \in L^2 _{\mathrm{loc}}(\mu_t;\R^d)$.
\item
There exists $\psi \in BV (\Omega \times [0,T)) \cap C ^{\frac12} ([0,T); L^1 (\Omega))$
such that the following hold.
\begin{enumerate}[label=\textnormal{(b\arabic*)}]
\item $\psi ^{\eps} \to \psi $ in $L^1  (\Omega \times [0,T))$
and a.e. pointwise, where $\psi ^{\eps}=\frac12 (\varphi^{\eps} +1)$.
\item $\psi \vert_{t=0}=\mathbbm{1} _{U_0}$ a.e. on $\Omega$.
\item For any $t \in [0,T )$, $\psi (\cdot ,t) =1$ or $0$ a.e. on $\Omega$ and
$\psi$ satisfies the volume preserving property, that is, 
\[
\int _{\Omega} \psi (x,t) \, \dx  = \mathcal{L}^{d} (U_0)
\qquad \text{for all} \ t \in [0,T).
\]
\item\label{(b4)} For any $t \in [0,T)$ and for any nonnegative $\phi \in C_c (\Omega\times [0,T) ;[0,\infty))$, 
we have $\| \nabla \psi (\cdot ,t) \| (\phi (\cdot,t) ) \leq \mu _t (\\phi (\cdot,t))$.
\end{enumerate}
\item\label{(c)} For $\lambda ^{\eps}$ given by \eqref{lambdadef}, we have
\begin{equation}\label{bound L^2 lambda^eps}
    \sup_{\eps\in(0,\eps_1)} \int _{t_1} ^{t_2} \vert \lambda ^{\eps} \vert^2 \, \dt <c(1+t_2-t_1)
\end{equation}
for every $0\leq t_1<t_2\leq T$, for some $\eps_1=\eps_1(d,D_0,D_1,\alpha)>0$ and $c=c(d,D_0,D_1)>0$.
Moreover, there exists $\lambda \in L ^2 (0,T)$ such that
$\lambda ^{\eps} \to \lambda$ weakly in $L^2 (0,T)$.
\item\label{(d)} 
There exists $\vec{f} \in L ^2 ([0,T); (L^2 (\mu _t;\R^d)))$ such that
\begin{equation}\label{claim-e2}
\begin{split}
&\lim_{\eps \rightarrow 0}
\frac{1}{\sigma} \int_{t_1} ^{t_2} \int _{\Omega }
-\lambda ^{\eps} \sqrt{2W (\varphi ^{\eps})}
\nabla \varphi ^{\eps} \cdot \vec{\phi}  \, \dx\dt \\
= & \int_{t_1} ^{t_2} \int _{\Omega} \vec{f} \cdot \vec{\phi}  \, {\rm d} \mu_t \dt 
= 
\int_{t_1} ^{t_2} \int _{\Omega} -\lambda \vec{\nu}\cdot \vec{\phi}
\, {\rm d}  \| \nabla \psi (\cdot,t) \| \dt
\end{split}
\end{equation}
for any $\vec{\phi} \in C_c (\Omega \times [0,T); \R^d)$ and every $0\leq t_1<t_2\leq T$,
where $\vec \nu$ is the inner unit normal vector of $\{ \psi (\cdot,t) =1 \}$
 on $\operatorname{spt} \| \nabla \psi (\cdot,t) \|$.

\item\label{(e)} The family of Radon measures $\{ \mu _t \}_{t \in [0,T)}$
is an $L^2$-flow in the sense of Definition \ref{def L^2 flow} with a generalized velocity vector $\vec{v} = \vec{h} +\vec{f}$,
where $\vec{h} \in L ^2 (0,T ; L ^2 (\mu_t;\R^d) ) $ 
is the generalized mean curvature vector of $\mu _t$.
Moreover, for every $\vec{\phi}  \in C_c (\Omega \times [0,T); \R^d)$ and
$\phi \in C_c (\Omega \times [0,T); [0,\infty))$,
\begin{equation}\label{claim-e}
\lim_{\eps\rightarrow 0}\frac{1}{\sigma}\int_{0}^{T}\int_\Omega  \vec{v^\eps} \cdot\vec{\phi} \,\eps|\nabla\varphi^\eps|^2\dx\dt
= \int_{0} ^T \int _{\Omega} \vec{v} \cdot \vec{\phi}  \, {\rm d} \mu_t \dt,
\end{equation}
and, for almost every $t\in (0,T)$,
\begin{equation}\label{lim h^eps}
- \int _{\Omega} \vec{\phi}  \cdot \hh (\cdot,t) \, {\rm d}\mu _t
=\frac{1}{\sigma}\lim _{\eps\rightarrow 0}
\int _{\Omega} \eps(\vec{\phi} \cdot \nabla \varphi ^{\varepsilon})
h^\eps
\dx,
\end{equation}
\begin{equation}\label{liminf h^2}
\int _\Omega \phi \vert \vec{h} \vert^2 \, {\rm d}\mu _t
\leq
\frac{1}{\sigma} \liminf_{\eps \rightarrow 0}
\int _{\Omega} \varepsilon  \phi 
|h^\eps|^2
\dx <\infty.
\end{equation}
\end{enumerate}

\end{theorem}

\begin{remark}
From \ref{(b4)}, \ref{(d)}, and \ref{(e)}, we have \eqref{v=h-lambda distrib}
for any $\vec{\phi} \in C_c (\Omega \times (0,T); \R^d)$, and hence $\vec{v} = \vec{h} -\lambda \frac{{\rm d} \| \nabla \psi (\cdot,t) \|}{\dmu _t} \vec{\nu}$ distributionally. 
\end{remark}

\begin{remark}\label{measurability}
The measurability of $\int_\Omega \vec\phi \cdot \vec h(\cdot,t) \, \dd \mu_t$ and $\int_\Omega \phi |\vec h(\cdot,t)|^2 \, \dd \mu_t$ with respect to time is non-trivial from its construction, regardless of the volume-preserving condition. 
Since the results of \cite{TakasaoVPMCF} implicitly assumed these measurabilities, we provide a proof for them in the Appendix using the argument from \cite{LiuWorkman}.
Below, we briefly explain why this measurability is non-trivial. Set $a_i (t):=\int_{\Omega} \varepsilon _{i} |h ^{\varepsilon_{i}} (\cdot,t)|^2 \, \dmu_t ^{\varepsilon_{i}}$ and $S(t):= \{ \{ \varepsilon_{i_j}\} \subset \{ \varepsilon_{i}\} \mid 
\sup_{j\geq 1} a_{i_j} (t)<\infty\}$.
By \cite[equation (62)]{TakasaoVPMCF} and Fatou's lemma,
\[
\int_0 ^T \liminf_{i\to\infty} a_i(t) \, \dt
\leq 
\liminf_{i\to\infty} \int_0 ^T a_i(t) \, \dt
\leq
\sup_{i\geq 1}\int_0 ^T a_i (t) \, \dt <\infty.    
\]
Hence for a.e. $t \in [0,T]$ we have
$\liminf_{i\to\infty} a_i(t) <\infty$ and hence $S(t)\not =\emptyset$. For such $t$,  
there exists a unique $(d{-}1)$-integral varifold $V_t\in \V_{d-1}(\Omega)$ such that $\mu_t=\|V_t\|$ and having generalized mean curvature vector $\vec{h} (\cdot,t) \in L^2 (\mu_t;\R^d)$ (for example, see \cite{TakasaoTonegawaMCFhigherdim}). 
If $\{\varepsilon_{i_j}\}, \{\tilde \varepsilon_{i_j}\} \in S(t)$, 
by the uniqueness of the limit, 
\begin{equation}\label{varifold-unique}
    \delta V_t ^{\varepsilon_{i_j}} (\vec g) \to \delta V_t (\vec g) \quad \text{and} \quad
\delta V_t ^{\tilde \varepsilon_{i_j}} (\vec g) \to \delta V_t (\vec g) \qquad\text{for any} \ \vec g \in C_c (\Omega).
\end{equation}
However, when $\{\varepsilon_{i_j}\} \notin S (t)$, it is non-trivial whether \eqref{varifold-unique} holds.
If there exists a measure zero set $N$ such that $\cap _{t \in [0,T] \setminus N} S(t) \not=\emptyset$, then there exists $\{\varepsilon_{i_j}\} \in \cap _{t \in [0,T] \setminus N} S(t)$ such that $\delta V_t ^{\varepsilon_{i_j}} (\vec g) \to \delta V_t (\vec g)$ for a.e. $t \in [0,T]$ and thus 
$\delta V_t (\vec g) = - \int \vec h (\cdot,t) \cdot \vec g \, \dmu_t$ is a measurable with respect to $t$ (since it is the limit of the measurable functions $f_j (t):=\delta V_t ^{\varepsilon_{i_j}} (\vec g)$).
However, in general, $\liminf_{i\to \infty} a_i (t) <\infty$ for a.e. $t \in[0,T]$ does not imply $\cap _{t \in [0,T] \setminus N} S(t) \not=\emptyset$. 
\end{remark}

\subsection{Main result: Existence of a volume-preserving Brakke-flow}\label{sec: mai results Brakke}

\begin{theorem}[Existence of a weak Brakke-flow]\label{thm: existence Brakke}
    Under the assumptions of Theorem \ref{mainthm1}, there exists a constant $C>0$ (depending on $T, d, D_0$, and $D_1$) such that the family of Radon measures $\{ \mu _t \}_{t \in [0,T)}$ defined in Theorem \ref{mainthm1} satisfies \eqref{Brakke ineq weak} with the constant $C$. In particular, $\{\mu_t\}_{t\in[0,T)}$ is a volume-preserving Brakke-flow in the sense of Definition \ref{def Brakke flow}.
\end{theorem}

\begin{proof}
First, let us extend $\varphi^\eps$, $\mu^\eps_t$, and $\mu_t$ periodically on $\R^d$. By integration by parts and the definition of $h^\eps$ \eqref{approximate v and h}, we have (cf. \cite[(57)]{TakasaoVPMCF})
    \begin{align*}
        \sigma\frac{\rm d}{{\rm d}t}\mu^\eps_t(\phi)&=\int_{\R^d} -\eps\phi (\varphi^\eps_t)^2+\lambda^\eps\phi\sqrt{2W(\varphi^\eps)}\varphi^\eps_t-\eps(\nabla \phi{\cdot}\nabla\varphi^\eps)\varphi^\eps_t\dx+\sigma\int_{\R^d} \phi_t\dmu^\eps_t\\
        &\stackrel{\eqref{approximate v and h}}{=}\int_{\R^d} -\eps\phi \left(h^\eps+\lambda^\eps\frac{\sqrt{2W(\varphi^\eps)}}{\eps}\right)^2+\lambda^\eps\phi\sqrt{2W(\varphi^\eps)}\left(h^\eps+\lambda^\eps\frac{\sqrt{2W(\varphi^\eps)}}{\eps}\right)\dx\\
        &\quad-\int_{\R^d}\eps(\nabla \phi{\cdot}\nabla\varphi^\eps)\left(h^\eps+\lambda^\eps\frac{\sqrt{2W(\varphi^\eps)}}{\eps}\right)\dx+\sigma\int_\Omega \phi_t\dmu^\eps_t\\
        &=\int_{\R^d} -\eps\phi\left(|h^\eps|^2+h^\eps\lambda^\eps\frac{\sqrt{2W(\varphi^\eps)}}{\eps}\right)-\eps\nabla \phi{\cdot}\nabla\varphi^\eps\left(h^\eps+\lambda^\eps\frac{\sqrt{2W(\varphi^\eps)}}{\eps}\right)\dx\\
        &\quad+\sigma\int_{\R^d} \phi_t\dmu^\eps_t\\
        &=\int_{\R^d} -\eps\phi\left(\frac{|h^\eps|^2}{2}+\frac{1}{2}\left|h^\eps+\lambda^\eps\frac{\sqrt{2W(\varphi^\eps)}}{\eps}\right|^2\right)+ \phi|\lambda^\eps|^2\frac{W(\phi^\eps)}{\eps}\dx\\
        &\quad-\int_{\R^d}\eps\nabla \phi{\cdot}\nabla\varphi^\eps\left(h^\eps+\lambda^\eps\frac{\sqrt{2W(\varphi^\eps)}}{\eps}\right)\dx+\sigma\int_{\R^d} \phi_t\dmu^\eps_t\\
        &\stackrel{\eqref{approximate v and h}}{=}\int_{\R^d} -\eps\phi\left(\frac{|h^\eps|^2}{2}+\frac{1}{2}|\varphi^\eps_t|^2\right)+ \phi|\lambda^\eps|^2\frac{W(\phi^\eps)}{\eps}\dx\\
        &\quad-\int_{\R^d}\eps\nabla \phi{\cdot}\nabla\varphi^\eps\left(h^\eps+\lambda^\eps\frac{\sqrt{2W(\varphi^\eps)}}{\eps}\right)\dx+\sigma\int_{\R^d} \phi_t\dmu^\eps_t   
    \end{align*}
    for every $\phi\in C^1_c({\R^d}\times[0,T];[0,\infty))$. Integrating in time between $0\leq t_1< t_2\leq T$, we find
    \begin{align}
        \sigma \mu^\eps_{t}(\phi)\Big|_{t=t_1}^{t_2}&=\int_{t_1}^{t_2}\int_{\R^d} -\eps\phi\left(\frac{|h^\eps|^2}{2}+\frac{1}{2}|\varphi^\eps_t|^2\right)+ \phi|\lambda^\eps|^2\frac{W(\varphi^\eps)}{\eps}\dx\dt\notag\\
        &\quad-\int_{t_1}^{t_2}\int_{\R^d}\eps\nabla \phi{\cdot}\nabla\varphi^\eps\left(h^\eps+\lambda^\eps\frac{\sqrt{2W(\varphi^\eps)}}{\eps}\right)\dx\dt+\sigma\int_{t_1}^{t_2}\int_{\R^d} \phi_t\dmu^\eps_t\dt.\label{Brakke eps} 
    \end{align}
    Passing to the limit for $\eps\rightarrow 0$, the left-hand side and the last term of the right-hand side converge by Theorem \ref{mainthm1} \ref{(a)}. For the third term of the right-hand side by Theorem \ref{mainthm1} \ref{(b4)}, \ref{(d)}, and \ref{(e)}, we get
    \begin{align*}
        &\lim_{\eps\rightarrow 0}-\int_{t_1}^{t_2}\int_{\R^d}\eps\nabla \phi{\cdot}\nabla\varphi^\eps\left(h^\eps+\lambda^\eps\frac{\sqrt{2W(\varphi^\eps)}}{\eps}\right)\dx\dt\\
    &\quad\stackrel{\eqref{claim-e2},\eqref{lim h^eps}}{=}\sigma\int_{t_1}^{t_2}\int_{\R^d}\nabla \phi{\cdot}\left(\hh-\lambda\nu\frac{{\rm d}\|\nabla \psi(\cdot,t)\|}{\dmu_t}\right)\dmu_t\dt
    \end{align*}    
    For the remaining terms on the right-hand side of \eqref{Brakke eps}, we first notice that the $\lim$ is less or equal than the $\limsup$, which is less or equal than the sum of the $\limsup$. Let $\eta \in C_c((t_1,t_2);[0,\infty))$ with $0\leq \eta\leq 1$ and set
    \begin{equation*}
        \zeta(t)\coloneqq \eta(t)\int_{\R^d}\phi(x,t)\frac{|\vec{h}(x,t)|^2}{2}\dmu_t(x).
    \end{equation*}
    We may assume without loss of generality $\operatorname{spt}\phi(\cdot,t)\subset\Omega +\vec{u}$ for every $t\in (0,T)$, for some $\vec{u}\in \R^d$. Hence, by periodicity of $\mu_t$ and the fact that $\vec{h}\in L^2(0,T;L^2(\mu_t;\R^d))$, we have
    \begin{align*}
        |\zeta(t)|\leq \frac{\|\phi\|_{L^\infty}}{2}\int_{\Omega+\vec{u}}|\vec{h}(x,t)|^2\dmu_t(x)=\frac{\|\phi\|_{L^\infty}}{2}\int_{\Omega}|\vec{h}(x,t)|^2\dmu_t(x)\in L^1(0,T).
    \end{align*}
    Thus, by Theorem \ref{mainthm1} \ref{(e)} and the dominated convergence theorem, 
    \begin{align*}
        \limsup_{\eps\rightarrow 0}\int_{t_1}^{t_2}\int_{\R^d} -\eps\phi\frac{|h^\eps|^2}{2}\dx\dt&=-\liminf_{\eps\rightarrow0}\int_{t_1}^{t_2}\int_{\R^d} \eps\phi\frac{|h^\eps|^2}{2}\dx\dt\\
        &\leq -\liminf_{\eps\rightarrow0}\int_{0}^{T}\int_{\R^d} \eps\eta\phi\frac{|h^\eps|^2}{2}\dx\dt\stackrel{\eqref{liminf h^2}}{\leq} -\sigma\int_{0}^{T} \int_{\R^d} \eta\phi\frac{|\vec{h}|^2}{2}\dmu_t\dt,
    \end{align*}
    which by letting $\eta \rightarrow \mathbbm{1}_{[t_1,t_2]}$ yields the desired limit.
    For the velocity term let us define the measure $d\widetilde{\mu}^\eps\coloneqq \frac{\eps}{\sigma}|\nabla \varphi^\eps|^2dxdt = d\mu^\eps +d\xi^\eps$, where $\xi^\eps$ is the discrepancy measure, given by 
    \begin{equation*}
        \xi^\eps(\zeta)\coloneqq \frac{1}{\sigma}\int_0^T\int_{\R^d} \zeta \left(\frac{\eps|\nabla\varphi^\eps|^2}{2}-\frac{W(\varphi^\eps)}{\eps}\right)\dx\dt\quad \forall \, \zeta \in C_c (\R^d \times [0,\infty)).
    \end{equation*}
    By Theorem \ref{mainthm1} \ref{(a)}, and the fact that $\xi^\eps\rightarrow0$ as Radon measures, we have $\widetilde{\mu}^\eps\rightarrow\mu$  as Radon measures, where  ${\rm d}\mu={\rm d}\mu_t\dt$.  By definition of approximate velocity \eqref{approximate v and h} and of $\widetilde{\mu}^\eps$,  we have
    \begin{align*}
        \int_{0}^{T}\int_{\Omega} |v^\eps|^2{\rm d}\widetilde{\mu}^\eps\leq 
        \frac{1}{\sigma}\int_{0}^{T}\int_{\Omega}\eps|\varphi^\eps_t|^2\dx\leq \frac{c}{\sigma}\mu^\eps_0(\Omega)\leq \frac{c}{\sigma} D_1,
    \end{align*}
    where the second inequality follows by the standard energy estimate for the solution to \eqref{ac}, cf. \cite[Prop. 7]{TakasaoVPMCF}. Hence, by \cite[Thm. 4.4.2]{Hutchinson} (or see \cite[Thm. B.3]{MugnaiRoeger}), there exists $g\in L^2(\dmu;\R^d)$ such that (up to a nonrelabeled subsequence)
\begin{align*}
        \int_{0}^{T}\int_{\R^d} g\cdot \zeta \dmu_t\dt&=\lim_{\eps\rightarrow 0}\frac{1}{\sigma}\int_{0}^{T}\int_{\R^d}v^\eps\cdot\zeta \,\eps|\nabla\varphi^\eps|^2\dx\dt
        \qquad \forall \, \zeta \in C_c (\R^d\times (0,T);\R^d).
    \end{align*}
By Theorem \ref{mainthm1} \ref{(e)}, it follows that $g=v$. Moreover, again by \cite[Thm. 4.4.2]{Hutchinson},
    \begin{align*}
        \limsup_{\eps\rightarrow0}-\frac{1}{2\sigma}\int_{t_1}^{t_2}\int_{\R^d}\eps\phi|\varphi^\eps_t|^2\dx= -\liminf_{\eps\rightarrow0}\frac{1}{2\sigma}\int_{t_1}^{t_2}\int_{\R^d}\eps\phi|\varphi^\eps_t|^2\dx\leq -\frac{1}{2}\int_{t_1}^{t_2}\int_{\R^d}\phi|v|^2\dmu_t\dt.
    \end{align*}

    It thus only remains to bound the $\limsup$ of the $|\lambda^\eps|^2$-term in \eqref{Brakke eps}. Let $\operatorname{spt}\phi\subset\subset B_r(x_0)$ for some $r\in (0,1)$ and some $x_0\in {\R^d}$. We have 
    \begin{align}
        \frac{1}{\sigma}\int_{t_1}^{t_2}\int_{\R^d} \phi|\lambda^\eps|^2\frac{W(\varphi^\eps)}{\eps}\dx\dt&\leq \int_{t_1}^{t_2}|\lambda^\eps|^2\mu^\eps_t(\phi)\dt\leq \|\phi\|_{L^\infty}\int_{t_1}^{t_2}|\lambda^\eps|^2\mu^\eps_t(B_r(x_0))\dt.\label{bound L^2 norm lambda^eps }
    \end{align}
We recall that by \cite[Corollary 1]{TakasaoVPMCF}, we have the following upper estimate on the density of $\mu^\eps_t$
\begin{equation*}
    \mu^\eps_t(B_r(x_0))\leq cr^{d-1},
\end{equation*}
for every $x_0\in \R^d$, $t\in [0,T]$, and for small enough $\eps$.
 Hence, by the $L^2$-bound \eqref{bound L^2 lambda^eps} of Theorem \ref{mainthm1} \ref{(c)}  
\begin{align*}
        \frac{1}{\sigma}
        \int_{t_1}^{t_2}\int_{\R^d} \phi|\lambda^\eps|^2\frac{W(\varphi^\eps)}{\eps}\dx\dt\leq r^{d-1}C(1+t_2-t_1)\|\phi\|_{L^\infty}
\end{align*}
and we have concluded the proof.
\end{proof}

\appendix

\section{Measurability of the varifolds with respect to time}

Let $\{\mu_t\}_{t\in [0,\infty)}$ be the family of Radon measures obtained in \cite{TakasaoVPMCF}. The following argument is based on Liu-Workman~\cite{LiuWorkman}. 

Let $N_1 \subset [0,\infty)$ be a set such that 
\begin{enumerate}
\item $\mathcal{L}^1 (N_1)=0$,
\item for all $t \in [0,\infty)\setminus N_1$ there exists a unique integral varifold $V_t$ with $\|V_t\| = \mu_t$
\item for all $t \in [0,\infty)\setminus N_1$, there exists a generalized mean curvature vector $\vec{h} (\cdot,t)$ of $V_t$, that is, 
$\delta V_t (\vec\eta) = - \int \vec{h}(\cdot,t) \cdot \vec \eta \, \dd \mu_t $ for any $\vec \eta \in C(\mathbb{T}^d ; \mathbb{R}^d)$. In addition, 
$ \int |\vec{h} (\cdot,t) | ^2 \, \dd \mu_t <\infty $.
\end{enumerate}
From the discussion in \cite{TakasaoVPMCF}, it is easy to check that such an $N_1$ exists.



\begin{proposition}
\begin{enumerate}
\item \label{measurability1} There exists a countable set $N_2 \subset [0,\infty)$ such that
for any $\phi \in C (\mathbb{T}^d;[0,\infty))$, $t \mapsto \mu _t (\phi)$ is a continuous function on $[0,\infty) \setminus N_2$. 
\item \label{measurability2} For any $\phi \in C (\mathbb{T}^d)$, $t \mapsto \mu _t (\phi)$ is a measurable function.
\item \label{measurability4} For any $\eta \in \mathbb{G}_{d-1} (\mathbb{T}^d)$, 
$t \mapsto V_t (\eta)$ is a measurable function.
\item \label{measurability5} For any $\vec{\eta} \in C (\mathbb{T}^d)$, 
$t \mapsto \int \vec{h} (\cdot,t) \cdot \vec{\eta} \, \dd \mu_t$ is a measurable function.
\item \label{measurability6} For any $\phi \in C (\mathbb{T}^d; [0,\infty))$, 
$t \mapsto \int \phi |\vec{h}(\cdot,t)|^2 \, \dd \mu_t$ is a measurable function.
\item \label{measurability7} For any $\phi \in C (\mathbb{T}^d)$, 
$t \mapsto \int \phi |\vec{h}(\cdot,t)|^2 \, \dd \mu_t$ is a measurable function.
\end{enumerate}
\end{proposition}

\begin{proof}
(\ref{measurability1}) Let $\{ \phi_k \}_{k=1} ^\infty \subset C(\mathbb{T}^d;[0,\infty))$ be a dense set with $\phi_k \in C^2 (\mathbb{T}^d;[0,\infty))$ for any $k \geq 1$. Since $t \mapsto  \mu_t (\phi _k)$ can be written as the difference of two increasing functions, there exists a countable set $C_k$ such that $t \mapsto  \mu_t (\phi _k)$ is continuous on $[0,\infty) \setminus C_k$ (see \cite{TakasaoVPMCF}). Set $N_2 = \cup_{k=1} ^\infty C_k$. Then, for any $\phi \in C (\mathbb{T}^d;[0,\infty))$, a standard $\varepsilon$-$\delta$ argument shows that $t \mapsto  \mu_t (\phi)$ is continuous on $[0,\infty) \setminus N_2$.

\bigskip

\noindent
(\ref{measurability2}) This claim follows immediately from (\ref{measurability1}).

\bigskip

\noindent
(\ref{measurability4}) Set $A:= [0, \infty)\setminus (N_1 \cup N_2)$.
Fix an arbitrary $\eta \in \mathbb{G}_{d-1} (\mathbb{T}^d)$. 
For any $i \in \mathbb{N}$, we define 
\[
A_i := \left\{ t \in A \mid  \int |\vec{h} (\cdot,t)|^2 \, \dd \mu_t \leq i \right\}.
\]
Fix an arbitrary $r\in \R$ and denote 
\[
A_i (r):=\{ t \mid V_t (\eta) \geq r \} \cap A_i.
\]
We will prove that $A_i (r)$ is relatively closed in $A$. We only need to show that if $\{t_j \}_{j=1} ^\infty \subset A_i (r)$ and $s \in A$ satisfy $t_j \to s$, then $s \in A_i (r)$.
Assume that $\{t_j \}_{j=1} ^\infty \subset A_i (r)$ and $s \in A$ satisfy $t_j \to s$.

Since $\sup_{j \geq 1}\int |\vec{h} (\cdot,t_j)|^2 \, \dd \mu_t \leq i$, $\sup_{j\geq 1}\mu_{t_j} (\mathbb{T}^d) <\infty$, and there exists $c>0$ such that
$\Theta^{\ast (d-1)} (\mu_{t_j}, x) \geq c$ for $\mu_{t_j}$-a.e. $x$ and for any $j\geq 1$, the compactness theorem of the rectifiable varifolds \cite{Allard} yields that there exist a rectifiable varifold $\tilde V$ and a subsequence $t_{j_k} \to s$ such that $V_{t_{j_k}} \stackrel{\ast}{\rightharpoonup} \tilde V$. 
In addition, by $s \in A$, $\mu_t (\eta)$ is continuous at $t=s$ for any $\eta \in C(\mathbb{T}^d)$ and thus $\mu_{t_{i_j}}\stackrel{\ast}{\rightharpoonup} \mu_s$. Since 
$\|\tilde V\|=\mu_s = \| V_s \|$
and there is only one rectifiable varifold $\tilde V$ satisfying $\mu_s = \| \tilde V \|$, we have $\tilde V=V_s$.
Therefore $V_{t_{j_k}} \stackrel{\ast}{\rightharpoonup} V_s$ implies
\[
r\leq \lim_{k\to\infty} V_{t_{j_k}} (\eta) = V_s (\eta).
\]
In addition, by the lower semi-continuity of the measure-function pair \cite{Hutchinson}, we have
\[
\int |\vec{h} (\cdot,s)|^2 \, \dd \mu_s \leq
\liminf_{k\to \infty} \int |\vec{h} (\cdot,t_{j_k})|^2 \, \dd \mu_{t_{j_k}} \leq i.
\]
Hence $s \in A_i (r)$. 
Since $A_i (r)$ is relatively closed in $A$, there exists a closed set $B_i (r) \subset [0,\infty)$ such that $A_i (r) = A \cap B_i (r)$, and thus $A_i (r)$ is measurable(because $A$ and $B_i (r)$ are measurable). Since
\[
\{ t \mid V_t (\eta) \geq r \} \cap A = \cup_{i=1} ^\infty A_i (r), \qquad \mathcal{L}^1 ([0, \infty)\setminus A)=0,
\]
$\{ t \mid V_t (\eta) \geq r \}$ is measurable. Hence $t \mapsto V_t (\eta)$ is a measurable function.

\bigskip

\noindent
(\ref{measurability5}) Fix an arbitrary $\vec{\eta} \in C^1 (\mathbb{T}^d)$.
We compute
\[
\int \vec{h} (\cdot,t) \cdot \vec{\eta} \, \dd \mu_t = - \int \nabla \vec{\eta} (x) \cdot S \, \dd V_t (x,S) = V_t (\phi),
\]
where $\phi (x,S)= \nabla \vec \eta (x) \cdot S$.
Hence by (\ref{measurability4}), $t \mapsto \int \vec{h}(\cdot,t) \cdot \vec{\eta} \, \dd \mu_t$ is a measurable function.
Fix an arbitrary $\vec{\eta} \in C (\mathbb{T}^d)$.
By a smooth approximation $\vec\eta _k$ of $ \vec \eta$, one can show that 
$t \mapsto \int \vec{h}(\cdot,t) \cdot \vec{\eta} \, \dd \mu_t$ is a measurable function.

\bigskip

\noindent
(\ref{measurability6})
Fix an arbitrary $\phi \in C (\mathbb{T}^d;[0,\infty))$ and set $\phi_\varepsilon:=\phi +\varepsilon$ for $\varepsilon>0$. Set
\[
f_\varepsilon (t):=\int \phi_\varepsilon |\vec{h} (\cdot,t)|^2 \, \dd \mu_t.
\]
Since $\mathcal{L}^1 ([0,\infty)\cap A^c)=\mathcal{L}^1 (N_1 \cup N_2)=0$, we may restrict the function $f_\varepsilon (t)$ to $A$. 
We will prove 
\begin{equation}\label{appendix-b2}
f_\varepsilon (s)\leq
\liminf_{t \to s, \ t\in A}f_\varepsilon (t) \qquad \text{for any} \ s \in A.
\end{equation}
Let $\{ t_j \}_{j=1} ^\infty \subset A$ be a sequence such that $t_j \to s \in A$ as $j\to \infty$
and $\liminf_{t \to s, \ t\in A}f_\varepsilon (t) =\lim_{j\to \infty}f_\varepsilon (t_j)$. We may assume that there exist $C>0$ and $c>0$ such that $ \lim_{j\to \infty} f_\varepsilon (t_j) <C$, and $\Theta^{\ast (d-1)} (\mu_{t_i},x) \geq c$ for $\mu_{t_i}$-a.e. $x$, for any $i\geq 1$. By an argument similar to the above, the compactness theorem of the rectifiable varifolds \cite{Allard} and the uniqueness of the rectifiable varifolds yield that there exists a subsequence $t_{j_k} \to s$ such that $V_{t_{j_k}} \stackrel{\ast}{\rightharpoonup} V_s$.
Then, by the lower semi-continuity of the measure-function pair, we have 
\begin{equation*}
f_\varepsilon (s)
=\int \phi_\varepsilon |\vec{h} (\cdot,s)|^2 \, \dd \mu_s
\leq
\liminf_{k\to \infty} \int \phi_\varepsilon |\vec{h} (\cdot,t_{j_k})|^2 \, \dd \mu_{t_{j_k}}
\leq
\lim_{i \to \infty}f_\varepsilon (t_i)
=\liminf_{t \to s, \ t\in A}f_\varepsilon (t).  
\end{equation*}
Hence we obtain \eqref{appendix-b2}.
Since $f_\varepsilon (t)$ is a lower semi-continuous function on $A$ and thus is a measurable function. Therefore $f (t) = \lim_{\varepsilon \to 0} f_\varepsilon (t)$ is a measurable function.

\bigskip

\noindent
(\ref{measurability7}) This claim follows immediately from (\ref{measurability6}).
\end{proof}

\section*{Acknowledgments}

The authors would like to thank the anonymous referee for their careful reading of the manuscript and valuable suggestions.
AC is supported by the Austrian Science Fund (FWF) projects 10.55776/F65,  10.55776/I5149,
10.55776/P32788,     and   10.55776/I4354,
as well as by the OeAD-WTZ project CZ 09/2023. 
KT is supported by JSPS KAKENHI Grant Numbers
JP23K03180, JP23H00085, and
JP24K00531.
Part of this research was funded by the Mobility Fellowship of the University of Vienna and conducted during a visit to Kyoto University, whose warm hospitality is gratefully acknowledged. We thank Prof. Yoshihiro Tonegawa for sharing his thoughts on these topics and pointing out relevant references.
Remark \ref{measurability} was incorporated thanks to a valuable comment by Yutaka Hirayama, for which the authors express their sincere gratitude. We would like to thank Katerina Nik for drawing our attention to the paper \cite{LiuWorkman}.


\begin{thebibliography}{99}

\bibitem{Allard}
W. K.~Allard.
On the first variation of a varifold.
\emph{Ann. of Math. (2)} 95 (1972),
417--491.


\bibitem{Almgren}
F.~Almgren, J.~Taylor, L.~Wang.
Curvature-driven flows: a variational approach.
\emph{SIAM J. Control Optim.}   31 (1993), no. 2, 387--438.

\bibitem{Ambrosio}
L.~Ambrosio, N.~Fusco, D.~Pallara.
Functions of bounded variation and free discontinuity problems.
\emph{The Clarendon Press, Oxford University Press, New York.} (2000).

\bibitem{BMR}
M.~Bonforte, F.~Maggi, D.~Restrepo.
Asymptotic behavior of a diffused interface volume-preserving mean curvature flow.
Arch. Ration. Mech. Anal., 250 (2026), Paper No. 23, 65.

\bibitem{Brakke}
K.~A.~Brakke.
\emph{The motion of a surface by its mean curvature}.
Math. Notes, 20
Princeton University Press, Princeton, NJ, 1978. i+252 pp.


\bibitem{Ecker}
K.~Ecker. 
\emph{Regularity theory for mean curvature flow.}
Progr. Nonlinear Differential Equations Appl., 57
Birkhäuser Boston, Inc., Boston, MA, 2004.

\bibitem{EscherSimonett}
J.~Escher, G.~Simonett.
The volume preserving mean curvature flow near spheres.
\emph{Proc. Amer. Math. Soc.}   126 (1998), no. 9, 2789--2796.

\bibitem{Fischer}
J. Fischer, S. Hensel, T. Laux and T. M. Simon. 
The local structure of the energy landscape in multiphase mean curvature flow: Weak-strong uniqueness and stability of evolutions. 
preprint (2020), https://arxiv.org/abs/2003.05478. 

\bibitem{Gage}
M.~Gage.
\emph{On an area-preserving evolution equation for plane curves. Nonlinear problems in geometry} (Mobile, Ala., 1985), 51–62.
Contemp. Math., 51
American Mathematical Society, Providence, RI, 1986. 

\bibitem{Golovaty}
D.~Golovaty.
The volume-preserving motion by mean curvature as an asymptotic limit of reaction-diffusion equations.
\emph{Quart. Appl. Math.}   55 (1997), no. 2, 243--298.

\bibitem{Huisken}
G.~Huisken. 
The volume preserving mean curvature flow.
\emph{J. Reine Angew. Math.}   382 (1987), 35--48.

\bibitem{Hutchinson}
J. E.~Hutchinson.
Second fundamental form for varifolds and the existence of surfaces minimising curvature.
\emph{Indiana Univ. Math. J.} 35 (1986), 45--71.


\bibitem{Julin-Niinikoski}
V.~Julin, J.~Niinikoski.
Quantitative {A}lexandrov theorem and asymptotic behavior of the volume preserving mean curvature flow.
\emph{Anal. PDE} 
16 (2023), 679--710.


\bibitem{KasaiTonegawa}
K.~Kasai, Y.~Tonegawa.
A general regularity theory for weak mean curvature flow.
\emph{Calc. Var. Partial Differential Equations}
50 (2014), No. 1-2,
1--68.

\bibitem{KimKwon}
I.~Kim, D.~Kwon.
Volume preserving mean curvature flow for star-shaped sets.
\emph{Calc. Var. Partial Differential Equations}   59 (2020), no. 2, 40 pp.

\bibitem{LM}
B.~Lambert, E.~M\"{a}der-Baumdicker.
Nonlocal estimates for the volume preserving mean curvature flow and applications.
\emph{Calc. Var. Partial Differential Equations} 
62 (2023), Paper No. 202, 32.

\bibitem{LauxWeakStrongUniqueness}
T.~Laux.
Weak-strong uniqueness for volume-preserving mean curvature flow.
\emph{Rev. Mat. Iberoam.}   40 (2024), no. 1, 93--110.

\bibitem{LauxOtto}
T.~Laux, F.~Otto.
Convergence of the thresholding scheme for multi-phase mean-curvature flow.
\emph{Calc. Var. Partial Differential Equations} 55 (2016), no. 5, Art. 129, 74 pp.

\bibitem{LauxVolumePreserving}
T.~Laux, D.~Swartz.
Convergence of thresholding schemes incorporating bulk effects.
\emph{Interfaces Free Bound.} 19 (2017), no.2, 273–-304.

\bibitem{LauxSimon}
T.~Laux, T.~Simon.
Convergence of the Allen-Cahn equation to multiphase mean curvature flow.
\emph{Comm. Pure Appl. Math.}   71 (2018), no. 8, 1597–-1647.

\bibitem{LiuWorkman}
Y.~T.~Liu, M.~Workman. The space-time-Grassmann measure of the Brakke flow. preprint (2025), https://arxiv.org/abs/2512.19227.


\bibitem{LuckhausSturzenhecker}
S.~Luckhaus, T.~Sturzenhecker.
Implicit time discretization for the mean curvature flow equation.
\emph{Calc. Var. Partial Differential Equations}   3 (1995), no. 2, 253–271.

\bibitem{MugnaiRoeger}
L.~Mugnai, M.~Röger. 
The Allen-Cahn action functional in higher dimensions.
\emph{Interfaces Free Bound.} 10 (2008), no.1, 45--78.

\bibitem{MugnaiSeisSpadaro}
L.~Mugnai, C.~Seis, E.~Spadaro.
Global solutions to the volume-preserving mean-curvature flow.
\emph{Calc. Var. Partial Differential Equations}   55 (2016), no. 1, Art. 18, 23 pp.

\bibitem{Simon}
L.~Simon.
Lectures on geometric measure theory.
\emph{Proceedings of the Centre for Mathematical Analysis, Australian National University.} 3 (1983).

\bibitem{StuvardTonegawa}
S.~Stuvard, Y.~Tonegawa. 
On the existence of canonical multi-phase Brakke flows.
\emph{Adv. Calc. Var.}   17 (2024), no. 1, 33--78.

\bibitem{TakasaoVPMCF}
K.~Takasao.
The existence of a weak solution to volume preserving mean curvature flow in higher dimensions. 
\emph{Arch. Ration. Mech. Anal.} 247 (2023), no.3, Paper No. 52, 53 pp.

\bibitem{TakasaoTonegawaMCFhigherdim}
K.~Takasao, Y.~Tonegawa.
Existence and regularity of mean curvature flow with transport term in higher dimensions.
\emph{Math. Ann.} 364 (2016), no.3-4, 857--935.

\bibitem{Tonegawa2014}
Y.~Tonegawa.
A second derivative {H}\"older estimate for weak mean curvature flow.
\emph{Adv. Calc. Var.} 7 (2014) No. 1, 91--138.

\bibitem{TonegawaBook}
Y.~Tonegawa.
\emph{Brakke's mean curvature flow.
An introduction}.
Springer Briefs Math.
Springer, Singapore, 2019. xii+100 pp.



\end{thebibliography}
\end{document}